\documentclass[12pt]{article}

\usepackage[margin=1in]{geometry}  % set the margins to 1in on all sides
\usepackage{graphicx}              % to include figures
\usepackage{amsmath}               % great math stuff
\usepackage{amsfonts}              % for blackboard bold, etc
\usepackage{amsthm}                % better theorem environments

\newtheorem*{thma}{Theorem A}
\newtheorem*{thmb}{Theorem B}
\newtheorem*{thmc}{Theorem C}
\newtheorem*{thmd}{Theorem D}
\newtheorem*{thme}{Theorem E}
\newtheorem*{thmf}{Theorem F}
\newtheorem*{thm1}{Theorem 1}
\newtheorem*{thm2}{Theorem 2}
\newtheorem*{cor1}{Corollary 1}
\newtheorem*{cor2}{Corollary 2}
\newtheorem*{cor3}{Corollary 3}
\newtheorem*{lem}{Lemma}
\theoremstyle{definition}

\begin{document}

%\pagestyle{myheadings}
%\markboth{Eze R Nwaeze}{On the number of zeros of a polynomial in a specified disk}
\nocite{*}

\title{\bf On the number of zeros of a polynomial in a specified disk\footnote{This is a preprint of a paper whose final and definite form is published open access in the International Journal of Pure and Applied Mathematics.}}

\author{Eze R. Nwaeze\\
Department of Mathematics\\
Tuskegee University \\
Tuskegee, AL 36088, USA\\
e-mail: enwaeze@mytu.tuskegee.edu}
\date{}
\maketitle

\begin{abstract}
Let $p(z)=a_0+a_1z+a_2z^2+a_3z^3+\cdots+a_nz^n$  be a polynomial of degree $n,$ where the coefficients $a_j,$ $j \in \{0,1,2,\cdots n\},$ may be complex. We impose some restriction on the coefficients of the real part of the given polynomial and then estimate the maximum number of zeros such polynomial can possibly have in a specified disk.
\end{abstract}



\section{Introduction}
The classical Enestr\"om-Kakeya Theorem states the following
\begin{thma}
Let $p(z)=\displaystyle{\sum_{j=0}^{n}a_jz^j}$ be a polynomial with real coefficients satisfying $$0<a_0\leq a_1\leq a_2\leq a_3\ldots\leq a_n.$$ Then all the zeros of $p(z)$ lie in $|z|\leq 1.$
\end{thma}
By putting a restriction on the coefficients of a polynomial similar to that of the  Enestr\"om-Kakeya Theorem, Mohammad \cite{Moh} proved the following on the number of zeros that can be found in a specified disk.

\begin{thmb}\label{thm1.2}
Let $p(z)=\displaystyle{\sum_{j=0}^{n}a_jz^j}$ be a polynomial with real coefficients satisfying $0<a_0\leq a_1\leq a_2\leq a_3\ldots\leq a_n.$ Then the number of zeros  of $p$ in $|z|\leq \dfrac{1}{2}$ does not exceed $$1+\dfrac{1}{\log 2}\log\Big(\dfrac{a_n}{a_0}\Big)$$
\end{thmb}

In her dissertation work, Dewan \cite{Dew}  weakens the hypotheses of Theorem B and proved the following two results for polynomials with complex coefficients.

\begin{thmc}\label{thm1.3}
Let $p(z)=\displaystyle{\sum_{j=0}^{n}a_jz^j}$ be a polynomial such that $|\arg a_j-\beta|\leq\alpha\leq\dfrac{\pi}{2}$ for $j\in \{0,1,2,\ldots,n\}$ and for some real numbers $\alpha$ and $\beta,$ and $$0<|a_0|\leq |a_1|\leq |a_2|\leq |a_3|\ldots\leq |a_n|.$$ Then the number of zeros of $p$ in $|z|\leq 1/2$ does not exceed $$\dfrac{1}{\log 2}\log\dfrac{|a_n|(\cos\alpha+\sin\alpha+1)+2\sin\alpha\sum_{j=0}^{n-1}|a_j|}{|a_0|}.$$
\end{thmc}

\begin{thmd}\label{thm1.4}
Let $p(z)=\displaystyle{\sum_{j=0}^{n}a_jz^j}$  where $Re(a_j)=\alpha_j$ and $Im(a_j)=\beta_j$ for all $j$ and $0<\alpha_0\leq \alpha_1\leq \alpha_2\leq\cdots\leq \alpha_n.$ Then the number of zeros of $p$ in $|z|\leq1/2$ does not exceed $$1+\dfrac{1}{\log 2}\log\dfrac{\alpha_n+\sum_{j=0}^{n}|\beta_j|}{|a_0|}.$$
\end{thmd}
Pukhta  \cite{Puk} generalized Theorems C and D by finding the number of zeros in $|z|\leq\delta$ for $0<\delta<1.$ The next Theorem, due to Pukhta, deals with a monotonicity condition on the moduli of the coefficients.

\begin{thme}\label{thm1.5}
Let $p(z)=\displaystyle{\sum_{j=0}^{n}a_jz^j}$ be a polynomial such that $|\arg a_j-\beta|\leq\alpha\leq\dfrac{\pi}{2}$ for $j\in \{0,1,2,\ldots,n\}$ and for some real $\alpha$ and $\beta,$ and $$0<|a_0|\leq |a_1|\leq |a_2|\leq |a_3|\ldots\leq |a_n|.$$ Then the number of zeros of $p$ in $|z|\leq \delta,$ $0<\delta<1,$ does not exceed $$\dfrac{1}{\log 1/\delta}\log\dfrac{|a_n|(\cos\alpha+\sin\alpha+1)+2\sin\alpha\sum_{j=0}^{n-1}|a_j|}{|a_0|}.$$
\end{thme}
Pukhta \cite{Puk} also gave a result which involved a monotonicity condition on the real part of the coefficients. Though the proof presented by Pukhta is correct, there was a slight typographical error in the statement of the result as it appeared in print. The correct statement of the theorem is as follows.

\begin{thmf}\label{thm1.6}
Let $p(z)=\displaystyle{\sum_{j=0}^{n}a_jz^j}$  where $Re(a_j)=\alpha_j$ and $Im(a_j)=\beta_j$ for all $j$ and $0<\alpha_0\leq \alpha_1\leq \alpha_2\leq\cdots\leq \alpha_n.$ Then the number of zeros of $p$ in $|z|\leq\delta,$ $0<\delta<1,$ does not exceed  $$\dfrac{1}{\log 1/\delta}\log2\Bigg[\dfrac{\alpha_n+\sum_{j=0}^{n}|\beta_j|}{|a_0|}\Bigg].$$
\end{thmf}

In this paper we generalize Theorem F and prove the following.

\begin{thm1}\label{thm1.7}
Let $p(z)=\displaystyle{\sum_{j=0}^{n}a_jz^j},$ $a_0\neq 0,$ be a polynomial of degree $n$ with complex coefficients, $Re(a_j)=\alpha_j$ and $Im(a_j)=\beta_j$ for all $j.$ If for some real numbers $t,$ and for some $\lambda \in\{0,1,2,\cdots n\},$ $$ t+\alpha_{n}\leq \alpha_{n-1}\leq \ldots\leq\alpha_{\lambda}\geq\alpha_{\lambda - 1}\geq\alpha_{\lambda - 2}\geq \ldots\geq \alpha_{1}\geq \alpha_{0},$$
then the number of zeros of $p$ in $|z|\leq\delta,$ $0<\delta<1,$ does not exceed $$\dfrac{1}{\log 1/\delta}\log\dfrac{M_1}{|a_0|},$$
where $$M_1=|\alpha_0|-\alpha_0+|\alpha_n|-\alpha_n+|t|-t+2\alpha_{\lambda}+2\displaystyle{\sum_{j=0}^{n}|\beta_j|}.$$
\end{thm1}
 For $t=0$ we get the following.
\begin{cor1}\label{cor1}
Let $p(z)=\displaystyle{\sum_{j=0}^{n}a_jz^j},$ $a_0\neq 0,$ be a polynomial of degree $n$ with complex coefficients,  $Re(a_j)=\alpha_j$ and $Im(a_j)=\beta_j$ for all $j.$ If for some $\lambda \in\{0,1,2,\cdots n\},$ $$\alpha_{n}\leq \alpha_{n-1}\leq \ldots\leq\alpha_{\lambda}\geq\alpha_{\lambda - 1}\geq\alpha_{\lambda - 2}\geq \ldots\geq \alpha_{1}\geq \alpha_{0},$$
then the number of zeros of $p$ in $|z|\leq\delta,$ $0<\delta<1,$ does not exceed $$\dfrac{1}{\log 1/\delta}\log\dfrac{M_1}{|a_0|},$$
where $$M_1=|\alpha_0|-\alpha_0+|\alpha_n|-\alpha_n+2\alpha_{\lambda}+2\displaystyle{\sum_{j=0}^{n}|\beta_j|}.$$

\end{cor1}

If $\lambda =0,$ then Corollary \ref{cor1} reduces to

\begin{cor2}\label{cor2}
Let $p(z)=\displaystyle{\sum_{j=0}^{n}a_jz^j},$ $a_0\neq 0,$ be a polynomial of degree $n$ with complex coefficients, $Re(a_j)=\alpha_j$ and $Im(a_j)=\beta_j$ for all $j.$ Suppose $$\alpha_{n}\leq \alpha_{n-1}\leq \ldots\leq\alpha_{0},$$
then the number of zeros of $p$ in $|z|\leq\delta,$ $0<\delta<1,$ does not exceed $$\dfrac{1}{\log 1/\delta}\log\dfrac{M_1}{|a_0|},$$
where $$M_1=|\alpha_0|+\alpha_0+|\alpha_n|-\alpha_n+2\displaystyle{\sum_{j=0}^{n}|\beta_j|}.$$

\end{cor2}

If, also, $\lambda =n,$ then Corollary \ref{cor1} becomes

\begin{cor3}\label{cor3}
Let $p(z)=\displaystyle{\sum_{j=0}^{n}a_jz^j},$ $a_0\neq 0,$ be a polynomial of degree $n$ with complex coefficients, $Re(a_j)=\alpha_j$ and $Im(a_j)=\beta_j$ for all $j.$ Suppose $$\alpha_{n}\geq\alpha_{n - 1}\geq\alpha_{n - 2}\geq \ldots\geq \alpha_{1}\geq \alpha_{0},$$
then the number of zeros of $p$ in $|z|\leq\delta,$ $0<\delta<1,$ does not exceed $$\dfrac{1}{\log 1/\delta}\log\dfrac{M_1}{|a_0|},$$
where $$M_1=|\alpha_0|-\alpha_0+|\alpha_n|+\alpha_n+2\displaystyle{\sum_{j=0}^{n}|\beta_j|}.$$

\end{cor3}

%\begin{cor}\label{cor1.11}
%Let $p(z)=\displaystyle{\sum_{j=0}^{n}a_jz^j},$ $a_0\neq 0,$ be a polynomial of degree $n$ with complex coefficients. If $Re(a_j)=\alpha_j$ and $Im(a_j)=\beta_j$ for all $j.$ Suppose $$\alpha_{n}\geq\alpha_{n - 1}\geq\alpha_{n - 2}\geq \ldots\geq \alpha_{1}\geq \alpha_{0}>0,$$
%then the number of zeros of $p$ in $|z|\leq\delta,$ $0<\delta<1,$ does not exceed $$\dfrac{1}{\log 1/\delta}\log\dfrac{M_1}{|a_0|}$$
%where $$M_1=$$
%
%\end{cor}

Suppose we assume $\alpha_{0}>0,$ then  Corollary 3  becomes Theorem F. Instead of proving Theorem 1, we prove the following more general result.
\begin{thm2}
Let $p(z)=\displaystyle{\sum_{j=0}^{n}a_jz^j},$ $a_0\neq 0,$ be a polynomial of degree $n$ with complex coefficients, $Re(a_j)=\alpha_j$ and $Im(a_j)=\beta_j$ for all $j.$ If for some real numbers $t,$ $s,$ and for some $\lambda \in\{0,1,2,\cdots n\},$ $$ t+\alpha_{n}\leq \alpha_{n-1}\leq \ldots\leq\alpha_{\lambda}\geq\alpha_{\lambda - 1}\geq\alpha_{\lambda - 2}\geq \ldots\geq \alpha_{1}\geq \alpha_{0}-s,$$
then the number of zeros of $p$ in $|z|\leq\delta,$ $0<\delta<1,$ does not exceed $$\dfrac{1}{\log 1/\delta}\log\dfrac{M_2}{|a_0|},$$
where $$M_2=|\alpha_0|-\alpha_0+|\alpha_n|-\alpha_n+|t|-t+|s|+s+2\alpha_{\lambda}+2\displaystyle{\sum_{j=0}^{n}|\beta_j|}.$$
\end{thm2}

Clearly $M_2$ is nonnegative.

\section{Lemma}
For the proof of our result we shall make use of the following result (see page 171 of the second edition) \cite{Tit}.

\begin{lem}
Let $F(z)$ be analytic in $|z|\leq R.$ Let $|F(z)|\leq M$ in the disk $|z|\leq R$ and suppose $F(0)\neq 0.$ Then for $0<\delta<1$ the number of zeros of $F(z)$ in the disk $|z|\leq \delta R$ is less than $$\dfrac{1}{\log 1/\delta}\log\dfrac{M}{|F(0)|}.$$
\end{lem}

\section{Proof of the Theorem}
\begin{proof}
Consider the polynomial
\begin{align*}
g(z)&=(1-z)p(z)\\
&=-a_{n}z^{n+1}+\displaystyle{\sum_{j=1}^{n}(a_j-a_{j-1})z^j}+a_0
\end{align*}
For $|z|=1,$
\begin{align*}
|g(z)|& \leq |a_{n}|+\displaystyle{\sum_{j=1}^{n}|a_j-a_{j-1}|}+|a_0|\\
&\leq |\alpha_n|+|\beta_n|+\displaystyle{\sum_{j=1}^{n}|\alpha_j-\alpha_{j-1}|}+\displaystyle{\sum_{j=1}^{n}|\beta_j-\beta_{j-1}|}+|\alpha_0|+|\beta_0|\\
&\leq |\alpha_n|+|\alpha_0|+\displaystyle{\sum_{j=1}^{n}|\alpha_j-\alpha_{j-1}|}+2\displaystyle{\sum_{j=0}^{n}|\beta_j|}\\
&\leq |\alpha_n|+|\alpha_0|+\displaystyle{\sum_{j=2}^{n-2}|\alpha_j-\alpha_{j-1}|}+|\alpha_{n-1}-\alpha_{n-2}|+|\alpha_{n}-\alpha_{n-1}|+|\alpha_{1}-\alpha_{0}|+2\displaystyle{\sum_{j=0}^{n}|\beta_j|}\\
&\leq |\alpha_0|-\alpha_0+|\alpha_n|-\alpha_n+\alpha_{n-2}+\alpha_1+|t|-t+|s|+s+\displaystyle{\sum_{j=2}^{\lambda}|\alpha_j-\alpha_{j-1}|}+\displaystyle{\sum_{j=\lambda+1}^{n-2}|\alpha_j-\alpha_{j-1}|}\\
&+ 2\displaystyle{\sum_{j=0}^{n}|\beta_j|}\\
&= |\alpha_0|-\alpha_0+|\alpha_n|-\alpha_n+|t|-t+|s|+s+2\alpha_{\lambda}+2\displaystyle{\sum_{j=0}^{n}|\beta_j|}\\
&=M_2.
\end{align*}
 Now $g(z)$ is analytic in $|z|\leq 1,$ and $|g(z)|\leq M_2$ for $|z|=1.$ So by the above lemma and the Maximum Modulus Principle, the number of zeros of $g$ (and hence of $p$) in $|z|\leq \delta$ is less than or equal to $$\dfrac{1}{\log 1/\delta}\log\dfrac{M_2}{|a_0|}.$$
Hence, the theorem follows.
\end{proof}

{\bf Acknowledgement:} The author is thankful to the anonymous referee for his/her valuable suggestions.

\end{document}